# AUTOMORPHISMS OF THE SEMIGROUP OF ENDOMORPHISMS OF FREE ALGEBRAS OF HOMOGENEOUS VARIETIES

R. LIPYANSKI

ABSTRACT. We consider homogeneous varieties of linear algebras over an associative-commutative ring $K$ with 1, i.e., the varieties in which free algebras are graded. Let $F = F(x_1, ..., x_n)$ be a free algebra of some variety $\Theta$ of linear algebras over $K$ freely generated by a set $X = \{x_1, ..., x_n\}$, $End\, F$ be the semigroup of endomorphisms of $F$, and $Aut\, End\, F$ be the group of automorphisms of the semigroup $End\, F$. We investigate structure of the group $Aut\, End\, F$ and its relation to the algebraical and categorical equivalence of algebras from $\Theta$.

We define a wide class of $R_1MF$-domains containing, in particular, Bezout domains, unique factorization domains, and some other domains. We show that every automorphism $\Phi$ of semigroup $End\, F$, where $F$ is a free finitely generated Lie algebra over an $R_1MF$-domain, is semi-inner. This solves the Problem 5.1 left open in [21]. As a corollary, semi-innerity of all automorphism of the category of free Lie algebras over $R_1MF$-domains is obtained. Relations between categorical and geometrical equivalence of Lie algebras over $R_1MF$-domains are clarified.

The group $Aut\, End\, F$ for the variety of $m$-nilpotent associative algebras over $R_1MF$-domains is described. As a consequence, a complete description of the group of automorphisms of the full matrix semigroup of $n \times n$ matrices over $R_1MF$-domains is obtained.

We give an example of the variety $\Theta$ of linear algebras over a Dedekind domain such that not all automorphisms of $Aut\, End\, F$ are quasi-inner.

The results obtained generalize the previous studies of various special cases of varieties of linear algebras over infinite fields.

## 1. INTRODUCTION

Let us recall the main definitions from Universal Algebraic Geometry [23, 24]. Let $\Theta$ be a variety of algebras over a commutative-associative ring $K$ and $F = F(X)$ be a free algebra from $\Theta$ generated by a finite set $X$. Here $X$ is supposed to be a subset of some infinite universum $X^0$. The set $Hom(F, G)$, $G \in \Theta$, can be treated as an affine space whose points are homomorphisms. The algebraic set in $Hom(F, G)$ and the category $K_\Theta(G)$ of algebraic sets over $G$ can be defined. The category $K_\Theta(G)$ is a geometric invariant of algebra $G$. Algebras $G_1$ and $G_2$ from $\Theta$ are categorically equivalent if the categories $K_\Theta(G_1)$ and $K_\Theta(G_2)$ are correctly isomorphic. Algebras $G_1$ and $G_2$ are geometrically equivalent if

$$T''_{G_1} = T''_{G_2}$$

holds for all finite sets $X$ and for all binary relations $T$ on $F$ and $'$ is Galois correspondence between sets in $Hom(F, G)$ and the binary relations on $F$.

It has been shown in [24] that categorical and geometrical equivalences of algebras are related and their relation is determined by the structure of the group $Aut\,\Theta^0$, where $\Theta^0$ is the category of free finitely generated algebras of $\Theta$. Note that the category $\Theta^0$ is small. The group $Aut\,\Theta^0$ is known for the following varieties: the variety of all groups, the variety of $F$-groups, where $F$ is a free group of constants, the variety of all semigroups, the variety of commutative-associative algebras with unit element over infinite fields, the variety of associative algebras over infinite fields, the variety of all Lie algebras over infinite fields, the variety of modules over $IBN$-rings [2, 3, 15, 16, 20, 21, 27].

There is a natural connection between a structure of the groups $Aut\,End\,F$, $F \in \Theta$, and $Aut\,\Theta^0$. However, a problem of description of the group $Aut\,End\,F$ is more complicated and was solved only for the following varieties: the variety of inverse semigroup, the variety of semigroups, the variety of groups, the variety of associative-commutative algebras over infinite fields, and the variety of Lie algebras over infinite fields [4, 8, 19, 20, 27].

We define a class of $R_1MF$-domains, containing Bezout domains, unique factorization domains, and some other domains. Namely, a domain $K$ is called $R_1MF$-domain if each $n \times m$ matrix $A$ over $K$ of rank 1 can be represented as a product of an $n \times 1$ matrix by an $1 \times m$ matrix over $K$. Here by the rank of matrix we understand its rank over the quotient field $\tilde{K}$ of $K$.

Our aim here is to describe the group $Aut\,End\,F$ and, as a consequence, to obtain a description of the group $Aut\,\Theta^0$ for the variety of Lie algebras and the variety of nilpotent associative algebras over $R_1MF$-domains.

The main theorems are as follows:

**Theorem A.** *Let $\mathcal{L}$ be the variety of Lie algebras over an $R_1MF$-domain $K$ and $F = F(x_1, .., x_n)$ be a finitely generated free Lie algebra of $\mathcal{L}$. Then any automorphism of the group $Aut\,End\,F$ is semi-inner.*

This Theorem solves Problem 5.1 in [21] formulated there for the variety of Lie algebras over fields. The description of the group $Aut\,End\,F(x_1, x_2)$ for the variety of Lie algebras over infinite fields has been given in [21]. The group $Aut\,End\,F(x_1, .., x_n)$, $n \geq 2$, for the variety of Lie algebras over infinite fields was described in [27].

**Theorem B.** Let $\mathcal{N}_m$ be the variety of $m$-nilpotent ($m \geq 2$) associative algebras over an $R_1MF$-domain $K$ and $F_m = F_m(X)$, $|X| < \infty$, be a finitely generated free algebra of the variety $\mathcal{N}_m$. The following three statements hold:

1. If either

    (a) $|K| = p^k$, $k \geq 1$, $p \neq 2$, and $p^k | m - 1$, and $m = 2r$, $r > 1$,

  or

    (b) $|K| = 2^k$, $k \geq 1$, (i.e. $p = 2$), and $2^k \mid m - 1$,

then the group $Aut\,End\,F_m$ is generated by semi-inner, $p$-semi-inner, mirror and $p$-mirror automorphisms of $End\,F_m$.

2. If $|K| = \infty$ or $|K| = p^k, k \geq 1$, and either

    (c) $p^k | m - 1$ and $m = 2r + 1$, $r \geq 0$,

  or

    (d) $p^k \nmid m - 1$ and $m \neq 2$,

then the group $Aut\,End\,F_m$ is generated by semi-inner and mirror automorphisms of $End\,F_m$.

3. If $K$ is any $R_1$MF-domain and $m = 2$, i.e., the multiplication in $F_m$ is trivial, then any automorphism of the group $Aut\,End\,F_m$ is semi-inner.

From part (3) of this theorem follows easily

**Corollary 1.1.** Let $M_n(K)$ be the full matrix semigroup of $n \times n$ matrices over a $R_1$MF-domain $K$. Then any automorphism of $M_n(K)$ is a semi-inner.

This result generalizes [13] on automorphisms of the full matrix semigroups of $n \times n$ matrices over principal ideal domains (see also [9], [10], where this result has been proved for the full matrix semigroup over fields).

Using Theorem A we further prove

**Theorem C.** Every automorphism of the category $\mathcal{L}^\circ$ of Lie algebras over an $R_1$MF-domain is semi-inner.

Earlier, this theorem has been proved for the case of Lie algebras over infinite fields [21]. Using Theorem B, a description of automorphisms of the category $\mathcal{N}_m^\circ$ can be derived.

We give also an example of variety $\Theta$ of linear algebras over a Dedekind domain for which the group $AutEnd\,F$ contains an automorphism which is not quasi-inner. Note that all automorphisms of $End\,F$ are quasi-inner in all above-mentioned varieties of algebras over $R_1$MF-domains.

The outline of this paper is as follows. We prove that all automorphisms of $End\,F$ for the varieties $\mathcal{B}_2$ and $\mathcal{N}_m$ are quasi-inner (see Theorem 5.3). Then we describe the bijections related to these quasi-inner automorphisms. To this end we investigate the structure of derivative algebras associated with quasi-inner automorphisms (see Propositions 6.6 and 6.7). Such a relation between quasi-inner automorphisms and derivative algebras was noted first in [27]. We prove the main statements of our paper, Theorems A, B, C, and discuss the relation between the categorical and geometrical equivalences of Lie algebras over $R_1$MF-domains (see Remark 7.1). Finally, we give the description of quasi-inner automorphisms of the semigroup $End\,A(X)$ for the variety of associative algebras over domains.

## 2. Automorphisms of the semigroup $End\,F$ and of the category $\Theta^0$

Throughout this paper "ring" will mean "commutative-associative ring with 1". Let $F = F(x_1, ..., x_n)$ be a free algebra of a variety $\Theta$ of linear algebras over ring $K$ generated by a set $X = \{x_1, ..., x_n\}$.

**Definition 2.1.** [2] An automorphism $\Phi$ of the semigroup $End\,F$ of endomorphisms of $F$ is called quasi-inner if there exists a bijection $s : F \to F$ such that $\Phi(\nu) = s\nu s^{-1}$, for any $\nu \in End\,F$; $s$ is called adjoint to $\Phi$.

**Definition 2.2.** [23] A quasi-inner automorphism $\Phi$ of $End\,F$ is called semi-inner if its adjoint bijection $s : F \to F$ satisfies the following conditions:
1. $s(a + b) = s(a) + s(b)$,
2. $s(a \cdot b) = s(a) \cdot s(b)$,
3. $s(\alpha a) = \varphi(\alpha)s(a)$,
for all $\alpha \in K$ and $a, b \in F$ and an automorphism $\varphi : K \to K$. If $\varphi$ is the identity automorphism of $K$, we say that $\Phi$ is an inner.

Let $A = A(x_1, ..., x_n)$ be a finitely generated free associative algebra of the variety $\mathcal{A}ss\text{-}\mathcal{K}$ of associative algebras over $K$.

**Definition 2.3.** [19] A quasi-inner automorphism $\Phi$ of $End\,A$ is called mirror if its adjoint bijection $s : A \to A$ is anti-automorphism of $A$.

Now we introduce a new class of quasi-inner automorphisms. Let $A = A(X)$ be a free finitely generated associative algebra over a ring $K$ of characteristic $p > 0$.

**Definition 2.4.** A quasi-inner automorphism $\Phi$ of $End\,A$ is called $p$-semi-inner ($p$-mirror) if $\Phi^p$ is a semi-inner (a mirror) automorphism of $End\,A$, whereas $\Phi^{p-1}$ is not.

Recall the notions of category isomorphism and equivalence [17]. An isomorphism $\varphi : \mathcal{C} \to \mathcal{D}$ of categories is a functor $\varphi$ from $\mathcal{C}$ to $\mathcal{D}$ which is a bijection both on objects and morphisms. In other words, there exists a functor $\psi : \mathcal{D} \to \mathcal{C}$ such that $\psi\varphi = 1_{\mathcal{C}}$ and $\varphi\psi = 1_{\mathcal{D}}$.

Let $\varphi_1$ and $\varphi_2$ be two functors from $\mathcal{C}_1$ to $\mathcal{C}_2$. A functor isomorphism $s : \varphi_1 \longrightarrow \varphi_2$ is a collection of isomorphisms $s_A : \varphi_1(A) \longrightarrow \varphi_2(A)$ defined for all $A \in Ob\,\mathcal{C}_1$ such that for every $\nu : A \longrightarrow B,\ \nu \in Mor\,\mathcal{C}_1,\ B \in Ob\,\mathcal{C}_1$, holds

$$s_B \cdot \varphi_1(\nu) = \varphi_2(\nu) \cdot s_A,$$

i.e., the following diagram is commutative

$$
\begin{array}{ccc}
\varphi_1(A) & \xrightarrow{\ s_A\ } & \varphi_2(A) \\
{\scriptstyle \varphi_1(\nu)}\downarrow & & \downarrow{\scriptstyle \varphi_2(\nu)} \\
\varphi_1(B) & \xrightarrow{\ s_B\ } & \varphi_2(B)
\end{array}
$$

The isomorphism of functors $\varphi_1$ and $\varphi_2$ is denoted by $\varphi_1 \cong \varphi_2$.

An equivalence between categories $\mathcal{C}$ and $\mathcal{D}$ is a pair of functors $\varphi : \mathcal{C} \to \mathcal{D}$ and $\psi : \mathcal{D} \to \mathcal{C}$ together with natural isomorphisms $\psi\varphi \cong 1_{\mathcal{C}}$ and $\varphi\psi \cong 1_{\mathcal{D}}$. If $\mathcal{C} = \mathcal{D}$, then we get the notions of automiorphism and autoequivalence of the category $\mathcal{C}$.

For every small category $\mathcal{C}$ denote the group of all its automorphisms by $Aut\,\mathcal{C}$. We will distinguish the following classes of automorphisms of $\mathcal{C}$.

**Definition 2.5.** [15, 21] An automorphism $\varphi : \mathcal{C} \to \mathcal{C}$ is equinumerous if $\varphi(A) \cong A$ for any object $A \in Ob\,\mathcal{C}$ ; $\varphi$ is stable if $\varphi(A) = A$ for any object $A \in Ob\,\mathcal{C}$ ; and $\varphi$ is inner if $\varphi$ and $1_{\mathcal{C}}$ are naturally isomorphic, i.e., $\varphi \cong 1_{\mathcal{C}}$.

In other words, an automorphism $\varphi$ is inner if for all $A \in Ob\,\mathcal{C}$ there exists an isomorphism $s_A : A \to \varphi(A)$ such that

$$\varphi(\nu) = s_B \nu s_A^{-1} : \varphi(A) \to \varphi(B)$$

for any morphism $\nu : A \to B$.

Let $\Theta$ be a variety of linear algebras over $K$. Denote by $\Theta^0$ the full subcategory of finitely generated free algebras $F(X), |X| < \infty$, of the variety $\Theta$.

**Definition 2.6.** [21] Let $A_1$ and $A_2$ be algebras from $\Theta$, $\delta$ be an automorphism of $K$ and $\varphi : A_1 \to A_2$ be a ring homomorphism of these algebras. A pair $(\delta, \varphi)$ is called semimomorphism from $A_1$ to $A_2$ if

$$\varphi(\alpha \cdot u) = \alpha^\delta \cdot \varphi(u), \quad \forall \alpha \in K,\ \forall u \in A_1.$$

Define the notion of a semi-inner automorphism of the category $\Theta^0$.

**Definition 2.7.** [21] An automorphism $\varphi \in Aut\,\Theta^0$ is called semi-inner if there exists a family of semi-isomorphisms $\{s_{F(X)} = (\delta, \tilde{\varphi}) : F(X) \rightarrow \tilde{\varphi}(F(X)),\ F(X) \in Ob\,\Theta^0\}$, where $\delta \in Aut\,K$ and $\tilde{\varphi}$ is a ring isomorphism from $F(X)$ to $\tilde{\varphi}(F(X))$ such that for any homomorphism $\nu : F(X) \longrightarrow F(Y)$ the following diagram

$$
\begin{array}{ccc}
F(X) & \xrightarrow{s_{F(X)}} & \tilde{\varphi}(F(X)) \\
\nu \downarrow & & \downarrow \varphi(\nu) \\
F(Y) & \xrightarrow{s_{F(Y)}} & \tilde{\varphi}(F(Y))
\end{array}
$$

is commutative.

Further, we will need the following

**Proposition 2.8.** [15, 21] For any equinumerous automorphism $\varphi \in Aut\,\mathcal{C}$ there exists a stable automorphism $\varphi_S$ and an inner automorphism $\varphi_I$ of the category $\mathcal{C}$ such that $\varphi = \varphi_S \varphi_I$.

## 3. Quasi-inner automorphisms of $End\,F$ for varieties of linear algebras

Now we introduce standard endomorphisms in free algebra $F = F(x_1, ..., x_n)$ of a variety $\Theta$.

**Definition 3.1.** Standard endomorphisms of $F$ in the base $X = \{x_1, ..., x_n\}$ are the endomorphisms $e_{ij}$ of $F$ which are determined on the free generators $x_k \in X$ by the rule: $e_{ij}(x_k) = \delta_{jk} x_i,\ x_i \in X,\ i, j, k \in [1n],\ \delta_{jk}$ is the Kronecker delta.

Denote by $S_0$ a subsemigroup of $End\,F$ generated by $e_{ij},\ i, j \in [1n]$.

**Proposition 3.2.** Let $\Phi \in Aut\,End\,F$. Elements of the semigroup $\Phi(S_0)$ are standard endomorphisms in some base $U = \{u_1, ..., u_n\}$ of $F$ if and only if $\Phi$ is a quasi-inner automorphism of $End\,F$.

*Proof.* Let $\Phi$ be a quasi-inner automorphism of $End\,F$. Consider the endomorphisms $\sigma$ and $\tau$ of $F$ given on generators $x_i \in X$ by the following rules: $\sigma(x_i) = s(x_i)$ and $\tau(x_i) = s^{-1}(x_i),\ i \in [1n]$. Let $\rho = \Phi(\tau)\sigma$. Then for any $x_i \in X$ we have:

$$\rho(x_i) = s\tau s^{-1}\sigma(x_i) = s\tau s^{-1}s(x_i) = s\tau(x_i) = x_i,\ i \in [1n],$$

i.e., $\rho = Id_F = \Phi(\tau)\sigma$, where $Id_F$ is the identical mapping on $F$. Replacing $\Phi$ by $\Phi^{-1}$ we obtain: $Id_F = \Phi^{-1}(\tau)\sigma$. Consequently, $Id_F = \tau\Phi(\sigma)$. Hence $\sigma$ is an automorphism of $F$.

Now we prove that $s(0) = 0$. For every $\zeta \in End\,F$ we have $\Phi(\zeta)(0) = s\zeta s^{-1}(0) = 0$, i.e., $\zeta s^{-1}(0) = s^{-1}(0)$. Assume that $\zeta$ is the zero endomorphism of $F$. Then $s^{-1}(0) = \zeta(s^{-1}(0)) = 0$.

Consider the following elements $u_1, ..., u_n$ from $F$:

$$u_1 = \sigma(x_1) = s(x_1), ..., u_n = \sigma(x_n) = s(x_n).$$

Since $\sigma$ is an automorphism of $F$, $u_1, ..., u_n$ is a base of $F$. Let us show that endomorphisms $\Phi(e_{ij}),\ i, j \in 1, .., n$ are standard endomorphisms in the base $U = \{u_1, ..., u_n\}$ :

$$\Phi(e_{ij})(u_k) = se_{ij}s^{-1}\sigma(x_k) = se_{ij}s^{-1}s(x_k) = se_{ij}(x_k) = s(\delta_{jk}x_i) = \delta_{jk}s(x_i) = \delta_{jk}u_i$$

Conversely, let $\Phi(e_{ij}),\ i, j \in [1n]$, be standard endomorphisms of $F$ in a base $U = \{u_1, ..., u_n\}$. Denote by $\mu_{ka},\ a \in F, k \in [1n]$, endomorphisms of $F$ given on generators $X$ by the rules: $\mu_{ka}(u_m) = \delta_{km}a$, $u_m \in U$. Then $\mu_{ka}e_{kk} = \mu_{ka}$. It is clear that if $\rho e_{kk} = \rho$ for some k and $\rho \in End\, F$, then there exists $a \in F$, such that $\rho = \mu_{ka}$.

Let, for $k = 1$, $\mu_{1a}e_{11} = \mu_{1a}$. Then $\Phi(\mu_{1a})\Phi(e_{11}) = \Phi(\mu_{1a})$. Since $\Phi(e_{11})$ is a matrix identity of $F$ in the base $U$, there exists an element $s(a) \in F$ such that $\Phi(\mu_{1a}) = \mu_{1s(a)}$. Note that $\Phi$ is an automorphism of $F$, hence $s$ is a bijection of $F$. Since $\rho\mu_{1a} = \mu_{1\rho(a)}$ for any $\rho \in End\, F$, we have $\Phi(\rho)\Phi(\mu_{1a}) = \Phi(\mu_{1\rho(a)})$. Therefore, $\Phi(\rho)\mu_{1s(a)}(u_1) = \mu_{1s(\rho(a))}(u_1)$. Thus, $\Phi(\rho)s(a) = s\rho(a), a \in F$, i.e. $\Phi(\rho) = s\rho s^{-1}$. $\qquad\square$

*Remark* 3.3. It is easy to show that the construction of the bijection $s : F \to F$ in the above proof does not depend on the choice of $k$.

*Remark* 3.4. From the proof of this Proposition we see also that the bijection $s$ transforms every base $X$ of $F$ into a base of the same algebra.

Now we define a notion of base $X_\sigma$-matrix of an automorphism $\Phi \in Aut End\, F$.

**Definition 3.5.** Let $\sigma$ be an element of the symmetric group $S_n$. The matrix $T_X^{(\sigma)} = (t_{ij}^{(\sigma)})$, where $t_{ij}^{(\sigma)} = \Phi(e_{ij})x_{\sigma(j)}$, $x_i \in X$, $i, j \in [1n]$, is called $X_\sigma$-matrix of $\Phi$ in base $X$.

The following Lemma establishes a useful property of $X_\sigma$-matrix of an automorphism $\Phi$ we need below.

**Lemma 3.6.** *Let $\alpha_1, ..., \alpha_n$ be endomorphisms of $F = F(x_1, ..., x_n)$. Then there exists an endomorphism $\alpha$ of $F$ such that*

$$(3.1) \qquad \alpha(t_{ij}^{(\sigma)}) = \alpha_i(t_{ij}^{(\sigma)})$$

*for all $i, j \in [1n]$ and some $\sigma \in S_n$. There exists a unique endomorphism $\alpha$ of $F$ such that (3.1) is fulfilled for all $i, j \in [1n]$ and all $\sigma \in S_n$.*

*Proof.* Let $\beta_i = \Phi^{-1}(\alpha_i)$ and $\beta_i(x_i) = y_i$, $i \in [1n]$. Determine an endomorphism $\beta \in End\, F$ on free generators $X$ of $F$ in the following way: $\beta(x_i) = y_i$, $i \in [1n]$ and let $\alpha = \Phi(\beta)$. Since $\beta e_{ij} = \beta_i e_{ij}$, we have $\Phi(\beta)\Phi(e_{ij})x_{\sigma(j)} = \Phi(\beta_i)\Phi(e_{ij})x_{\sigma(j)}$ for some $\sigma \in S_n$ and all $i, j \in [1n]$. Therefore, $\alpha(t_{ij}^{(\sigma)}) = \alpha_i(t_{ij}^{(\sigma)})$ for all $i, j \in [1n]$ and some $\sigma \in S_n$.

Let (3.1) be fulfilled for all $i, j \in [1n]$ and all $\sigma \in S_n$. We wish to prove the uniqueness of $\alpha$. Assume, on the contrary, that there exists $\gamma \in End\, F$ such that

$$\gamma(t_{ij}^{(\sigma)}) = \alpha_i(t_{ij}^{(\sigma)}), \ \text{ i.e., } \ \gamma\Phi(e_{ij})(x_{\sigma(j)}) = \alpha_i\Phi(e_{ij})(x_{\sigma(j)})$$

for all $i, j \in [1n]$ and all $\sigma \in S_n$. Then $\alpha\Phi(e_{ij}) = \gamma\Phi(e_{ij})$ for all $i, j \in [1n]$. Thus, $\Phi^{-1}(\alpha)e_{ij} = \Phi^{-1}(\gamma)e_{ij}$ and, as a consequence, $\Phi^{-1}(\alpha)e_{ij}x_j = \Phi^{-1}(\gamma)e_{ij}x_j$, i.e., $\Phi^{-1}(\alpha)x_i = \Phi^{-1}(\gamma)x_i$ for all $i$. We arrive at $\alpha = \gamma$. $\qquad\square$

## 4. $R_1MF$-domains

Let $K$ be an integral domain (domain, for short) and $\tilde{K} = Frac\, K$ be the quotient field of $K$.

**Definition 4.1.** Rank of a matrix $A \in M_{n \times m}(K)$ is the rank of $A$ over the field $\tilde{K}$, i.e., $rank_K\, A := rank_{\tilde{K}}\, A$.

**Definition 4.2.** We say that a domain $K$ satisfies the rank-1-matrix factorization condition ($R_1$MF-condition) if each $n \times m$ matrix $A$ over $K$ of rank 1 can be presented as a product of an $n \times 1$ matrix $c^i$ by an $1 \times m$ matrix $d_i$ over $K$, i.e., $A = c^i \cdot d_i$. A domain $K$ with $R_1$MF condition is called $R_1$MF-domain.

Now we will give several examples of $R_1$MF-domains.

**Example 4.3.** An $n \times m$ matrix $A$ is equivalent to an $n \times m$ matrix $B$ over a domain $K$ if there exist invertible matrices $P$ and $Q$ such that $A = PBQ$. Recall (see [14]) that an elementary divisor domain (EDD, for short) $K$ is a domain with the property that each matrix $A$ over $K$ is equivalent to a diagonal matrix

$$diag\,(d_1, d_2, ...) = \begin{pmatrix} d_1 & & & \\ & d_2 & & \\ & & \ddots & \\ & & & 0 \end{pmatrix},$$

where $d_i$ divides $d_{i+1}$ for all $i$.

For instance, any principal ideal domain is EDD. It is clear that each EDD is an $R_1$MF-domain.

**Example 4.4.** A Bezout domain is a domain in which any finitely generated ideal is principal, (see [6, 11]). By [6], Proposition 4.4, every Bezout domain is an $R_1$MF-domain. Note that every EDD is a Bezout domain. However, inverse inclusion is still an open question (see [11]).

**Example 4.5.** Let $K$ be a unique factorization domain (UFD, for short). Let us show that $K$ is an $R_1$MF-domain.

Let $A = (a_{ij}) \in M_{n \times m}(K)$ and $rank\,A = 1$. Then there exists $g \in K$ such that $gA = c^1 \cdot d_1$, where $c^1 = (c_{k1})$ and $d_1 = (d_{1k})$, $k \in [1n]$, $c_{k1}, d_{1k} \in K$. Thus $ga_{ij} = c_{i1}d_{1j}$. Assume that $g = p$ is a prime element in $K$. Since $K$ is UFD, we have $p|c_{i1}$ or $p|d_{1j}$. If $p|c_{i1}$ for all $i$ then our statement is true. Let there exist $s \leq n$ such that $p|c_{k1}$ for all $k < s$ but $p \nmid c_{s1}$. Then $p|d_{1j}$ for all $j$ and this yields the assertion.

Now, let $g$ be a non-prime element in $K$. We can represent $g = p_1...p_r$, $r > 1$, where all $p_i$ are prime elements in $K$. Using induction on $r$ we obtain the proof of this statement in the general case.

Note that the group algebra over field $P$ of the additive group of rational numbers (written multiplicatively) is a Bezout domain but not UFD (see [6], exercise 3.5 ). This suggests that the classes of UFDs and of $R_1$MF-domains do not coincide.

Now consider an example of domain which is not an $R_1$MF-domain.

**Example 4.6.** A domain $K$ with quotient field $\tilde{K}$ is called a Dedekind domain if it satisfies any of the following equivalent conditions:

(i) every ideal in $K$ is projective;

(ii) every nonzero ideal $C$ of $K$ is invertible (that is $CC^{-1} = K$, where $C^{-1} = \{x \in \tilde{K} | xa \subset K\}$.

As follows from Corollary 5.6 (or Remark 5.7) the Dedekind domain $K = \{a + b\sqrt{-5} \mid a, b \in \mathbf{Z}\}$ is not an $R_1$MF-domain.

It is proved in [14] that a domain $K$ is EDD if and only if every $2 \times 1$ and $2 \times 2$ matrices over $K$ are equivalent to a diagonal matrix. Our hypothesis is: a domain $K$ is $R_1$MF-domain if and only if each $2 \times 2$ matrix over $K$ of rank 1 can be represented as a column-by-row product.

Note that there is a notion of rank of matrix over any ring (see, for example, [6], [11]). Thus the notion of $R_1$MF-domain can be generalized to non-commutative rings. In future we are going to study such a generalization.

Let $V$ be a free module of finite rank $n$ over an $R_1$MF-domain K, $P$ be a subsemigroup of $End\,V$ generated by non-zero elements $P_{ij} \in End\,V$, $i, j \in [1n]$ such that $P_{ij}P_{mk} = \delta_{jm}P_{ik}$, $P_{ij} \neq 0$, for every $i, j, k, m \in [1n]$. Denote by $E$ a subsemigroup of $End\,V$ generated by the elementary matrices $E_{ij}$ in a basis $B = \langle v_1, ..., v_n \rangle$ of the free module $V$, i.e., $E_{ij}(v_k) = \delta_{jk}v_i$, $i, j, k \in [1n]$.

**Lemma 4.7.** *The semigroup $P$ and $E$ are conjugate via an automorphism $\rho$ of $V$.*

*Proof.* Denote by $V_{\tilde{K}} = V \bigotimes_K \tilde{K}$ a vector space over the quotient field $\tilde{K}$ of the domain $K$ Let $dim\,V_{\tilde{K}} = n$. Since $P_{11} \neq 0$, there exists $u \in V$ such that $u_1 = P_{11}u \neq 0$. Let $u_i = P_{11}u_1$, $i \in [1n]$. It is easy to check that the ordered set $B_1 = \langle u_i \in V | i \in [1n] \rangle$ forms a basis of the vector space $V_{\tilde{K}}$. Denote by $f_1(P_{ij})$, $i, j \in [1n]$ the matrices of the elements $P_{ij}$ in the basis $B_1$. Then $f_1(P_{ij}) = E_{ij}$, $i, j \in [1n]$. From this follows that the rank of every $f_1(P_{ij})$ over $K$ is 1. Since $K$ is an $R_1$MF-domain, we can represent $f_1(P_{ii}) = a^i b_i$, where $a^i = (a_{si})$, $a_{si} \in K$, $s \in [1n]$, is an $n \times 1$ matrix and $b_i = (b_{is})$, $b_{is} \in K$, $s \in [1n]$, is an $1 \times n$ matrix. Denote $A = (a_{ij})$ and $B = (b_{ij})$. Since

$$f_1(P_{ii})f_1(P_{jj}) = \delta_{ij}f_1(P_{ii}) = \sum_{k=1}^{n} b_{ik}a_{kj}a^i b_j,$$

we obtain $BA = I$. Hence, $A^{-1}f_1(P_{ii})A = E_{ii}, i \in [1n]$. We may assume that $f_1(P_{ii}) = E_{ii}$ for any $i$ in a basis $B_2 = \langle w_i \in V | i \in [1n] \rangle$ of the module $V$ over $K$. Since $P_{ij} = P_{ii}P_{ij}P_{jj}$, there exist $d_{ij} \in K$ such that $f_1(P_{ij}) = d_{ij}E_{ij}$. It is clear that $d_{ij}d_{mk} = \delta_{jm}d_{ik}$, $d_{ii} = 1$ for all $i, j, k, m \in [1n]$ and $d_{ij}$ are units in $K$. Consider $B_3 = \langle v_i \in V | v_i = d_{1i}^{-1}w_i, i \in [1n] \rangle$, a basis of the module $V$ over $K$. It is easy to check that the elements $P_{ij}$ can be represented by the elementary matrices in this basis. This completes the proof. □

## 5. Homogeneous varieties of algebras and quasi-inner automorphisms of $End\,F$

Let $\Theta$ be a variety of linear algebras over a ring $K$ and $T(\Theta)$ be its T-ideal.

**Definition 5.1.** [26] The variety $\Theta$ is called homogeneous if its T-ideal $T(\Theta)$ is a homogeneous.

It is known that the varieties of associative algebras, of nilpotent associative algebras, of Lie algebras, of alternative algebras, of Jordan algebras (if $1/2 \in K$) are homogeneous varieties ([26]).

Further, we consider only homogeneous varieties with free algebras without 1. A free algebra $F(X)$ of such a variety can be naturally decomposed as

$$(i)\ F(X) = \bigoplus_{k=1}^{\infty} F^{(k)} \text{ and } (ii)\ F^{(k)}F^{(m)} \subseteq F^{(k+m)},$$

where $F^{(k)}$ is a $K$-submodule generated by all monomial of $F(X)$ of degree $k$.

Denote by $F' = \bigoplus_{k=2}^{\infty} F^{(k)}$. It is easy to prove the following

**Lemma 5.2.** *If $\varphi \in End\, F$, then $\varphi(F') \subseteq F'$*

Let $\mathcal{B}_2$ be the variety of linear algebras over an $R_1MF$-domain defined by the identity $x_1^2 = 0$, and $\mathcal{N}_m$, $m \geq 2$, be the variety of nilpotent algebras of class $\leq m$ over an $R_1MF$-domain, i.e., for every algebra $G \in \mathcal{N}_m$ holds $G^m = 0$. We will say, somewhat not rigorously, that $\mathcal{N}_m$ is the variety of $m$-nilpotent algebras.

**Proposition 5.3.** Let $\Theta$ be a homogeneous variety of linear algebras over an $R_1MF$-domain $K$ and $F = F(x_1, ..., x_n)$ be its free algebra. If either $\Theta \subseteq \mathcal{B}_2$ or $\Theta \subseteq \mathcal{N}_m$ for some $m \geq 2$, then all automorphisms of the semigroup $End\, F$ are quasi-inner.

*Proof.* Let $\Phi$ be an automorphism of $End\, F$ and $\sigma \in S_n$. Consider the $X_\sigma$-matrix $T_X^{(\sigma)}$ of $\Phi$ in the base $X = \{x_1, ..., x_n\}$:

$$T_X^{(\sigma)} = (t_{ij}^{(\sigma)}),$$

where $t_{ij}^{(\sigma)} = \Phi(e_{ij})x_{\sigma(j)}$, $i, j \in [1n]$. Any element $t_{ij}^{(\sigma)}$ of $T_X^{(\sigma)}$ can be written in the following form:

(5.1)
$$t_{ij}^{(\sigma)} = m_{ij}^{(\sigma)} + g_{ij}^{(\sigma)},$$

where $m_{ij}^{(\sigma)} \in F^{(1)}$ is a linear part of $t_{ij}^{(\sigma)}$, and $g_{ij}^{(\sigma)} \in F' = \bigoplus_{k=2}^{\infty} F^{(k)}$.

Consider two cases

1. Let $\Theta \subseteq \mathcal{B}_2$. We shall show that there exists a non-zero element $m_{ij}^{(\sigma)}$ for some $i, j$ and $\sigma \in S_n$. Assume, on the contrary, that $m_{ij}^{(\sigma)} = 0$ for all $i, j$ and $\sigma \in S_n$. Consider the following $n$ endomorphisms of $F$: $\alpha_1 = e_{11}, ..., \alpha_n = e_{nn}$. Note that by our assumption $t_{ij}^{(\sigma)} = g_{ij}^{(\sigma)} \in F'$. As a consequence, $e_{ii}(t_{ij}^{(\sigma)}) = 0$. By Lemma 3.6, there exists a unique endomorphism $\alpha \in End\, F$ such that

$$\alpha(t_{ij}^{(\sigma)}) = \alpha_i(t_{ij}^{(\sigma)}) = e_{ii}(t_{ij}^{(\sigma)}) = 0 \,\, \forall i, j \in [1n], \,\forall \sigma \in S_n.$$

However, $e_{11}(t_{ij}^{(\sigma)}) = e_{22}(t_{ij}^{(\sigma)}) = 0$ for all $i, j \in [1n]$ and all $\sigma \in S_n$. We arrived at a contradiction with the uniqueness of the endomorphism $\alpha$.

Let us fix $i, j \in [1n]$ and $\sigma \in S_n$ for which $m_{ij}^{(\sigma)} \neq 0$ and write the elements $\Phi(e_{ij})x_{\sigma(k)}, k \in [1n]$, in the ordered base $X^{(\sigma)} = \langle x_{\sigma(1)}, ..., x_{\sigma(n)} \rangle$:

(5.2)
$$\Phi(e_{ij})x_{\sigma(k)} = a_{k1}^{(ij)}x_{\sigma(1)} + ... + a_{kn}^{(ij)}x_{\sigma(n)} + f_k^{(ij)}(x_1, ..., x_n),$$

where $a_{km}^{(ij)} \in K$ and $f_k^{(ij)}(x_1, ..., x_n) \in F'$. Denote by $M_n(K)$ the full matrix semigroup of $n \times n$ matrices over $K$. There exists a mapping $\psi : \Phi(S_0) \to M_n(K)$ such that $\psi(\Phi(e_{ij})) = A_{ij}^{(\sigma)}$, where $A_{ij}^{(\sigma)} = (a_{km}^{(ij)})$ is the matrix of the linear part of (5.2). By Lemma 5.2 and the equality $\Phi(e_{kl})\Phi(e_{ls}) = \Phi(e_{ks})$, we obtain that $\psi$ is a homomorphism from $\Phi(S_0)$ to $M_n(K)$. Thus, $A_{kl}^{(\sigma)}A_{ls}^{(\sigma)} = A_{ks}^{(\sigma)}$. Since $m_{ij}^{(\sigma)} \neq 0$, we have $A_{ij}^{(\sigma)} \neq 0$ for some $i, j \in [1n]$ and $\sigma \in S_n$. From this it follows that the matrices $A_{lm}^{(\sigma)} \neq 0$ for all $l, m \in [1n]$. Consequently, $\psi$ is a monomorphism. Denote $\Phi(S_0) = A$.

Let $V$ be a free module over $K$ with a basis $\tilde{X} = \langle x_1, ..., x_n \rangle$. By Lemma 4.7, the semigroup $A$ and $E$ are conjugate, i.e., there exists an automorphism $\rho$ of $V$ such that $\rho A_{ij}^{(\sigma)} \rho^{-1} = E_{ij}$, where $E_{ij}$, $i, j \in [1n]$ are elementary matrices over $K$ in a basis $Y = \langle y_1, ..., y_n \rangle$ and $\rho(x_i) = y_i$. Since $F$ is a Hopfian algebra, $Y$ is a base of $F$. The elements $\Phi(e_{i1})y_1, i \in [1n]$, can be represented in the base $Y$ as

(5.3) $\qquad \Phi(e_{i1})y_1 = y_i + g_i(y_1, ..., y_n), \; g_i \in F', \; i \in [1n].$

Denote $Z = \{z_i | z_i = \Phi(e_{i1})y_1, i \in [1n]\}$. Now we have to prove that the elements of $Z$ form a base of $F$. Let $e'_{ij}, i, j \in [1n]$ be the standard endomorphisms of $F$ in the base $Y$, i.e., $e'_{ij}(y_k) = \delta_{jk}y_i, i, j, k \in [1n]$. Applying $e'_{ii}$ to (5.3), we obtain

$$e'_{ii}(z_i) = y_i + \eta(y_i), \; \eta(y_i) \in F', \; i \in [1n].$$

Since $z_i = \Phi(e_{i1})y_1, \; i \in [1n]$, are elements of $Y_\varepsilon$-matrix ($\varepsilon$ is the identical substitution from $S_n$) of automorphism $\Phi$ in the base $Y$, by Lemma 3.6, there exists an endomorphism $\alpha \in End\, F$ such that

(5.4) $\qquad \alpha(z_i) = e'_{ii}(z_i) = y_i + \eta(y_i).$

Since $x_1^2 = 0$ is the identity in $\Theta$, we have $\eta(y_i) = 0$. By Lemma 3.6, there exists an endomorphism $\alpha \in End\, F$ such that

$$\alpha(z_i) = e'_{ii}(z_i) = y_i.$$

Since $F$ is a Hopfian algebra, the elements $z_i, i \in [1n]$, form a base of $F$.

2. Let $\Theta$ be a subvariety of variety $\mathcal{N}_m$. We shall show that there exists $m_{ij}^{(\sigma)} \neq 0$ for some $i, j$ and $\sigma \in S_n$. Assume that $m_{ij}^{(\sigma)} = 0$ for all $i, j$ and all $\sigma \in S_n$. We have

$$\Phi(e_{11})x_{\sigma(1)} = g_{11}^{(\sigma)}(x_1, ..., x_n),$$

where $g_{11}^{(\sigma)} \in F'$. Now we will prove, as an intermediate result, that $g_{11}^{(\sigma)} = 0$ for all $\sigma \in S_n$. Let $g_{11}^{(\sigma)} \neq 0$ for some $\sigma \in S_n$. Thus, $1 < deg\, g_{11}^{(\sigma)} \leq m - 1$. However, from the equality

$$\Phi(e_{11})x_{\sigma(1)} = \Phi(e_{11}^m)x_{\sigma(1)}$$

follows $deg\, g_{11}^{(\sigma)} > m$. This contradiction leads to $g_{11}^{(\sigma)} = 0$.

Finally, $\Phi(e_{11})x_{\sigma(1)} = m_{11}^{(\sigma)} + g_{11}^{(\sigma)} = 0$ for all $\sigma \in S_n$ and, as a consequence, $\Phi(e_{11}) = 0$, i.e., $e_{11} = 0$. Thus, we arrive at a contradiction again. Therefore, there exists $m_{ij}^{(\sigma)} \neq 0$ for some $i, j$ and $\sigma \in S_n$.

As in the case 1 we obtain the equality (5.4). From (5.4) follows

$$\alpha(z_i) = y_i \; (mod\, F'), \; \text{for all } i \in [1n].$$

It is well known (see [7, 18]) that every endomorphism $\tau$ of a finitely generated nilpotent algebra $G$ which induces an invertible linear transformation on the free $K$-module $G/G'$ is an automorphism of $G$. Taking this into account we obtain that $\alpha$ is an automorphism of $F$. Thus, we have $\alpha^{-1}(y_i) = z_i \; (mod\, F')$, i.e., the elements $z_i, i \in [1n]$ form a basis of $K$-module $F/F'$. Since $F$ is a finitely generated nilpotent algebra, the elements $z_i, i \in [1n]$, form a base of $F$ (see [18]).

Now in the both cases, $\Theta \subseteq \mathcal{B}_2$ or $\Theta \subseteq \mathcal{N}_m$, we have

$$\Phi(e_{ij})z_m = \Phi(e_{ij})\Phi(e_{m1})y_1 = \Phi(e_{ij}e_{m1})y_1 = \Phi(\delta_{jm}e_{i1})y_1 = \delta_{jm}z_i.$$

By Proposition 3.2, the automorphism $\Phi$ is quasi-inner as claimed.  $\square$

Now we consider the variety $\Upsilon$ of linear algebras with zero multiplication over a Dedekind domain such that the group $Aut\, End\, F$, where $F$ is a free two-generated algebra over $K$, contains a non-quasi-inner automorphism. This example is a modification of an example by Isaacs [12].

**Example 5.4.** Let $K = \{a + b\sqrt{-5} \mid a, b \in \mathbf{Z}\}$ be a Dedekind domain and $\Upsilon$ be a variety of linear algebras over $K$ with zero multiplication. Let $V = K^2$ be a free module over $K$. The module $V$ can be considered as a free two-generated algebra over $K$ of variety $\Upsilon$. Note that a semigroup $End\, V$ is the full matrix semigroup $M_{2\times 2}(K)$. We wish to construct an automorphism of the semigroup $M_{2\times 2}(K)$ which is not quasi-inner.

Let

$$m = \left( \begin{array}{cc} 1 + \sqrt{-5} & -2 \\ -2 & 1 - \sqrt{-5} \end{array} \right),$$

so that

$$m^{-1} = \frac{1}{2} \left( \begin{array}{cc} 1 - \sqrt{-5} & 2 \\ 2 & 1 + \sqrt{-5} \end{array} \right).$$

Note that $m^{-1} \notin M_{2\times 2}(K)$ but $m^{-1}xm \in M_{2\times 2}(K)$, i.e., $\Phi(x) = m^{-1}xm$, $x \in M_{2\times 2}(K)$ is an automorphism of the semigroup $M_{2\times 2}(K)$. We will show that this automorphism is not quasi-inner.

Assume that $\Phi$ is a quasi-inner automorphism of $End\, F$, i.e., there exists a bijection $s : V \to V$ on $V$ such that $\Phi(x) = s^{-1}xs$, $x \in M_{2\times 2}(K)$. Consequently, $s^{-1}xs = m^{-1}xm$ for all $x \in M_{2\times 2}(K)$. Then $\sigma x = x\sigma$ for all $x \in M_{2\times 2}(K)$, where $\sigma = ms^{-1}$ is a mapping from $V$ to $V$. Next we will prove that $\sigma = \alpha I$, $\alpha \in K$.

Consider the linear transformation $\gamma_a : V \to V$ defined by $\gamma_a(e_1) = a$, $\gamma_a(e_2) = 0$, where $e_1 = (1, 0)$, $e_2 = (0, 1)$ is a basis of the module $V$ and $a \in V$. Let $\sigma(e_1) = \alpha e_1 + \beta e_2$, $\alpha, \beta \in K$. Since $\sigma\gamma_a(e_1) = \gamma_a\sigma(e_1)$, we have $\sigma(a) = \gamma_a(\alpha e_1 + \beta e_2) = \alpha a$, i.e., $\sigma = \alpha I$, $\alpha \in K$.

Now we have $ms^{-1} = \sigma = \alpha I$ and thus $s^{-1} = \alpha m^{-1}$. Since $s^{-1}$ is the bijection on $V$, we obtain that $\det \alpha m^{-1}$ is a unit in $K$. However, only $\pm 1$ are units in $K$. Therefore $\det \alpha m^{-1} = \pm 1$, and from this, $\frac{1}{2}\alpha^2 = \pm 1$, i.e., $\alpha = \pm\sqrt{\pm 2}$. Since $\alpha \notin K$ we arrived at a contradiction, i.e., $\Phi$ is not quasi-inner.

Note that we have also $\Phi^2(x) = m^{-2}xm^2 = a^{-1}xa$, where

$$a = \left( \begin{array}{cc} \sqrt{-5} & -2 \\ -2 & -\sqrt{-5} \end{array} \right) \in M_{2\times 2}(K)$$

Since $\det a = 1$, the matrix $a^{-1}$ belongs to $M_{2\times 2}(K)$ as well. Therefore, $\Phi^2$ is a quasi-inner automorphism of $End\, F$, whereas $\Phi$ is not.

This example counts in favour of the following problem:

**Problem 5.5.** *Let $\Theta$ be a variety of linear algebras over a domain $K$ such that $\Theta \subseteq \mathcal{B}_2$ or $\Theta \subseteq \mathcal{N}_c$ for some $c$. Let $\Phi$ be an automorphism of the semigroup $End\, F$, where $F = F(x_1, ..., x_n)$ is a free $n$-generated algebra in $\Theta$. It is true that there exists a natural number $k(n)$ such that $\Phi^{k(n)}$ is a quasi-inner automorphism of $End\, F$?*

The example 5.4 leads also to the following statement:

**Corollary 5.6.** There exists a Dedekind domain which is not an $R_1MF$-domain.

*Proof.* Let the assumptions of the example 5.4 be fulfilled. Suppose that the Dedekind domain $K = \{a + b\sqrt{-5} \mid a, b \in \mathbf{Z}\}$ is an $R_1MF$-domain. Since the variety $\Upsilon$ belongs to $\mathcal{B}_2$, by Proposition 5.3 all automorphisms of $End\, F$, where $F$ is a finitely generated free algebra of $\Upsilon$, are quasi-inner. This fact contradicts to the example 5.4. $\square$

*Remark* 5.7. Note that Corollary 5.6 can be also proved by direct calculations.

Indeed, assume that $K = \{a + b\sqrt{-5} \mid a, b \in \mathbf{Z}\}$ is an $R_1$MF-domain. Take the matrix

$$c = \begin{pmatrix} 1 + \sqrt{-5} & 2 \\ 3 & 1 - \sqrt{-5} \end{pmatrix} \in M_{2\times 2}(K).$$

Since $rank\, c = 1$, we can represent the matrix $c$ as a column-by-row product:

$$\begin{pmatrix} 1 + \sqrt{-5} & 2 \\ 3 & 1 - \sqrt{-5} \end{pmatrix} = \begin{pmatrix} x_1 + y_1\sqrt{-5} \\ x_2 + y_2\sqrt{-5} \end{pmatrix} \cdot \begin{pmatrix} x_3 + y_3\sqrt{-5} & x_4 + y_4\sqrt{-5} \end{pmatrix},$$

where $x_i, y_i \in \mathbf{Z}$, $i = [1, 4]$. It can be shown that this system of polynomial equations of the second order has no solutions over $\mathbf{Z}$. We omit the calculations. This contradiction gives us a different proof of Corollary 5.6.

## 6. Derivative algebras associated with a quasi-inner automorphisms

Let $\Phi \in Aut\, End\, F$, $F \in \Theta$, be a quasi-inner automorphism of the semigroup $End\, F$ with the adjoint bijection $s : F \to F$. Our goal is to describe these bijections of $F$. For this purpose we reformulate universal algebra notations and results from [27] for the category of linear algebras.

By Remark 3.4, the bijection $s$ transforms every base $X = \{x_1, ..., x_n\}$ of $F$ into a base $Y = \{y_1, ..., y_n\}$ of this algebra: $s(x_i) = y_i$, $\forall i \in [1n]$. Consider an automorphism $\sigma : F \to F$ such that $\sigma(x_i) = y_i$, $\forall i \in [1n]$. Denote by $s_1 = \sigma s^{-1}$ a bijection of $F$. We have $s_1(x_i) = x_i$, $\forall i \in [1n]$. Define two automorphisms of $End\, F$:

$$\Phi_1(\nu) = s_1 \nu s_1^{-1} \text{ and } \Phi_2(\nu) = \sigma \nu \sigma^{-1}, \ \forall \nu \in End\, F.$$

Then $\Phi = \Phi_1^{-1}\Phi_2$. Therefore, it is sufficient to investigate quasi-inner automorphisms of $End\, F$ adjoint bijections of which preserve bases of $F$. It can be assumed that the adjoint bijection $s$ of $\Phi$ fixes the base elements $x_i \in X$, i.e., $s(x_i) = x_i$, $\forall i \in [1n]$.

Denote by $\theta_{a_1, ..., a_n}$, $a_i \in F$, an endomorphism of $F$ given on generators $X$ by the following rules:

$$\theta_{a_1, ..., a_n}(x_1) = a_1, ..., \theta_{a_1, ..., a_n}(x_n) = a_n.$$

Then

(6.1) $$\Phi(\theta_{a_1, ..., a_n}) = \theta_{s(a_1), ..., s(a_n)}.$$

From (6.1) follows that

(6.2) $$sF(x_{i_1}, ..., x_{i_k}) = F(x_{i_1}, ..., x_{i_k}), \ x_{i_s} \in X.$$

Denote by $F^* = \langle F; \circ, \bot, *, 0 \rangle$ a derivative algebra with the same support $F$ as the original algebra and with one nullary operation $0$ which coincides with $0$ of $F$, one unary operation $\circ$ and two binary operations, $\bot$ and $*$, determined in the following way:

1. $\alpha \circ a_1 = \theta_{a_1} s(\alpha \cdot x_1)$, $\forall \alpha \in K$, $\forall a_1 \in F$ and $x_1 \in X$,
2. $a_1 \bot a_2 = \theta_{a_1, a_2} s(x_1 + x_2)$, $\forall a_1, a_2 \in F$ and $x_1, x_2 \in X$,
3. $a_1 * a_2 = \theta_{a_1, a_2} s(x_1 \cdot x_2)$, $\forall a_1, a_2 \in F$ and $x_1, x_2 \in X$,

where the operations written on the right side of these formulas are the main operations in $F$. We shall say that the derivative algebra $F^*$ is associated with the automorphism $\Phi$.

Now, for completeness of presentation, we give the proof of the following statement.

**Proposition 6.1.** [27] The following hold.

1. Algebra $F^* \in \Theta$;
2. The bijection $s$ is an isomorphism of $F$ into $F^*$.

*Proof.* We will check the compatibility of the bijection $s$ with operation $*$. We have

$$\Phi(\theta_{a_1,a_2})s(x_1 \cdot x_2) = s\theta_{a_1,a_2}s^{-1}s(x_1 \cdot x_2) = s\theta_{a_1,a_2}(x_1 \cdot x_2) = s(a_1 \cdot a_2)$$

On the other hand by (6.1)

$$\Phi(\theta_{a_1,a_2})s(x_1 \cdot x_2) = \theta_{s(a_1),s(a_2)}s(x_1 \cdot x_2) = s(a_1) * s(a_2),$$

i.e., $s(a_1 \cdot a_2) = s(a_1) * s(a_2)$. In a similar way it is easy to check the compatibility of $s$ with operations $\perp$ and $\circ$. Finally, since $s$ is a bijection on $F$, we have $F^* \in \Theta$ and $s : F \to F^*$ is an isomorphism. $\square$

Now using these results we investigate the structure of bijections $s$ adjoint to quasi-inner automorphisms $\Phi$ of $End\, F$ for some classes of homogeneous varieties of linear algebras.

Let $\Theta$ be a homogeneous variety of linear algebras with or without 1 over a ring $K$ such that for each its free algebra $F(X)$ the following condition holds: there exists a free associative algebra $U(F)$ (an enveloping algebra) containing the algebra $F$ and freely generated by the same set $X$. We will denote by "$\cdot$" an operation of multiplication of elements from $U(F)$ (for brevity, we will omit this sign if it is clear from context). Note that the variety of Lie algebras, the variety generated by the free special Jordan algebras and, of course, variety of associative algebras over a ring $K$ are such varieties [1, 7, 26]. We shall call these varieties $\mathcal{A}$-varieties of linear algebras.

Let $\Psi$ be a homogeneous variety of $m$-nilpotent linear algebras over a ring $K$ such that for each its free algebra $F$ freely generated by $X$ there exists a nilpotent associative algebra $U_m(F)$ containing the algebra $F$ and freely generated by the same set $X$. We shall call these varieties $\mathcal{A}_m$-varieties of linear algebras. It is known [5] that the variety $\mathcal{N}_m$ of $m$-nilpotent Lie algebras is an $\mathcal{A}_m$-variety.

Let $A(X)$ be a free associative algebra freely generated by $X$ over a ring $K$. Take

$$f = \sum a_I x_I,$$

where $x_I = x_{i_1}...x_{i_m}$, $x_{i_s} \in X$, are monomials in $A$ and $a_I \in K$ are almost all 0. Note that the empty product of $x_i \in X$ to represent 1. By the support of $f$ we understand the set of all $x_I$ such that $a_I \neq 0$.

Let $A = A(x_1, x_2, x_3)$ be a 3-generated associative algebra over a ring $K$ and $f(x_1, x_2)$ be a monomial in $A$. Denote by $M_f$ the support of $f(x_1 + x_3, x_2)$ in $A$. Now we need the following

**Lemma 6.2.** *If* $f(x_1, x_2)$ *and* $g(x_1, x_2)$ *be two different monomials in* $A$, *then*

$$M_f \bigcap M_g = \varnothing.$$

*Proof.* Let

$$f(x_1, x_2) = x_1^{k_{11}} x_2^{k_{12}}...x_1^{k_{1s}} \text{ and } g(x_1, x_2) = x_1^{k_{21}} x_2^{k_{22}}...x_1^{k_{2s}},$$

where $s \geq 1$ and $k_{ij} \in \mathbf{N}_0 = \mathbf{N} \cup \{0\}$. Let $r$ be the minimal natural number such that $k_{1r} \neq k_{2r}$. We consider separately even and odd natural numbers $r$. If $r = 2k + 1, k \geq 0$, our statement follows from the fact that all monomials in polynomials $(x_1 + x_3)^{k_{1,2k+1}}$ and $(x_1 + x_3)^{k_{2,2k+1}}$ are different. If $r = 2k, k \geq 1$, all monomials in polynomials $(x_1 + x_3)^{k_{1,2k-1}} x_2^{k_{1,2k}}$ and $(x_1 + x_3)^{k_{1,2k-1}} x_2^{k_{2,2k}}$ are different and the statement also follows. The proof is complete. $\qquad \square$

*Remark* 6.3. Clearly, a similar assertion can be formulated for the support of the polynomial $f(x_1, x_2 + x_3)$.

**Definition 6.4.** A polynomial $P(x_1, x_2) \in F(x_1, ..., x_n)$, $n \geq 2$, is called a distributive polynomial over $K$ if

1. $P(a + b, c) = P(a, c) + P(b, c)$,
2. $P(a, b + c) = P(a, b) + P(a, c)$

for any $a, b, c \in F(x_1, ..., x_n)$.

Consider distributive polynomials in a free algebra $F(x_1, ..., x_n)$ over a ring $K$ for an $\mathcal{A}$-variety $\Theta$.

**Lemma 6.5.** *If $P(x_1, x_2)$ is a distributive polynomial over $K$ in an $\mathcal{A}$-variety $\Theta$, then*

$$P(x_1, x_2) = \alpha x_1 x_2 + \beta x_2 x_1, \ \alpha, \beta \in K,$$

*is a representation of $P(x_1, x_2)$ in $U_n = U(F(x_1, ..., x_n))$, $n \geq 2$, where $F = F(x_1, ..., x_n)$ is a free algebra of $\Theta$.*

*Proof.* It is clear that any distributive polynomial contains no constant term. Write $P(x_1, x_2) = \sum_i \alpha_i f_i(x_1, x_2)$, where $\alpha_i \in K$ and $f_i(x_1, x_2) \in U_n$ are different monomials in $U_n$. By definition 6.4 we have

$$(6.3) \qquad \begin{aligned} P(x_1 + x_3, x_2) &= P(x_1, x_2) + P(x_3, x_2), \\ P(x_1, x_2 + x_3) &= P(x_1, x_2) + P(x_1, x_3), \ x_i \in X. \end{aligned}$$

Since $U_n$ is a free associative algebra, we have by Lemma 6.2

$$(6.4) \qquad \begin{aligned} (a) \quad & f_i(x_1 + x_3, x_2) = f_i(x_1, x_2) + f_i(x_3, x_2), \\ (b) \quad & f_i(x_1, x_2 + x_3) = f_i(x_1, x_2) + f_i(x_1, x_3) \end{aligned}$$

for all $i$. Let for some $i$

$$(6.5) \qquad f_i(x_1, x_2) = x_1^{k_{i1}} x_2^{k_{i2}} x_1^{k_{i3}} ... x_1^{k_{i,s-1}} x_2^{k_{is}}, s \geq 1.$$

We will prove that, in the representation (6.5), every $k_{ij} \in \{0, 1\}$. Assume, on the contrary, that there exists $k_{i,2m-1} > 1$ for some $m \geq 1$. Then the monomial

$$x_1 \underbrace{x_3 x_3 \ldots x_3}_{k_{i1}-1} x_2^{k_{i2}} \ldots x_1 \underbrace{x_3 x_3 \ldots x_3}_{k_{i,2m-1}-1 \neq 0} x_2^{k_{i,2m}} \ldots x_2^{k_{is}}$$

appears in (6.4 (a)) on the left but does not on the right. This gives a contradiction.

If there exists $k_{i,2m} > 1$, $m \geq 1$, in (6.5), then using the equality (6.4 (b)) we achieve a contradiction in a similar way. Therefore, every $k_{ir}$, $r \in [1s]$, is equal to 0 or 1.

Assume that $k_{i1} = 1$. Then arguing as above, we obtain $k_{i,2k+1} = 0$ for all $k \geq 1$. If all $k_{i,2k} = 0$, $k \geq 1$, then $f_i = x_1$ and we arrive at a contradiction with (6.4 (b)). Thus, without loss of generality, it can be assumed that $k_{i2} = 1$. As above, we obtain $k_{i,2k} = 0$ for all $k > 1$. Therefore, $f_i(x_1, x_2) = x_1 x_2$. Now, if $k_{i1} = 0$, then, without loss of generality, we can assume that $k_{i2} = 1$. As above, we obtain $f_i(x_1, x_2) = x_2 x_1$. Finally, we have our assertion. $\qquad\square$

We are now ready to prove

**Proposition 6.6.** Let $\Theta$ be an $\mathcal{A}$-variety of algebras over a ring $K$ and $\Phi \in AutEnd\,F$ be a quasi-inner automorphism of a semigroup $End\,F$. Let $F^*$ be the derivative algebras associated with $\Phi$. Then there exist $\alpha, \beta \in K$, $\alpha^2 + \beta^2 \neq 0$, such that the following statements hold

(i)   $a * b = \alpha a \cdot b + \beta b \cdot a$;

(ii)  $a \perp b = a + b$;

(iii) $\xi \circ a = \varphi(\xi) a$,

for any $a, b \in F^*$ and $\xi \in K$ and an automorphism $\varphi : K \to K$.

*Proof.* We begin by showing (i). Let $s$ be a bijection adjoint to $\Phi$. Since $sF(x_1, x_2) = F(x_1, x_2)$ (see the equality (6.2)), $s(x_1 x_2)$ is a polynomial $P(x_1, x_2)$ which belongs to $F(x_1, x_2)$. By definition of the operation $*$ in $F^*$ we have

$$a * b = \theta_{a,b} s(x_1 \cdot x_2) = \theta_{a,b} P(x_1, x_2) = P(a, b)$$

Since $*$ is a distributive operation in $F^*$, the polynomial $P(x_1, x_2)$ is distributive. By Lemma 6.5

$$P(x_1, x_2) = \alpha x_1 x_2 + \beta x_2 x_1, \ \alpha, \beta \in K.$$

If $\alpha = \beta = 0$, the derivative algebra $F^*$ is an algebra where the multiplication is trivial. Since $F^* \in \Theta$ and $\Theta$ is an $\mathcal{A}$-variety, we achieve a contradiction and the result follows.

(ii) By definition of operation $\perp$, we have

$$(6.6) \qquad x_1 \perp x_2 = s(x_1 + x_2) = G(x_1, x_2),$$

where $G(x_1, x_2) \in F(x_1, x_2)$. Thus, $x_1 = s(x_1) = s(x_1 + 0) = G(x_1, 0)$ and, similarly, $x_2 = G(0, x_2)$. Now we write $G(x_1, x_2)$ as an element of $U_2 = U(F(x_1, x_2))$

$$(6.7) \qquad G(x_1, x_2) = x_1 + x_2 + \gamma g(x_1, x_2),$$

where $\gamma \in K$, $g(x_1, x_2) \in F(x_1, x_2) \subseteq U_2$ and the degree of $g(x_1, x_2)$ in $U_2$ is equal to $k \geq 2$. We will show that, in fact, $\gamma = 0$.

Assume, on the contrary, that $\gamma \neq 0$. Consider the equality

$$(6.8) \qquad x_1 * (x_2 \perp x_3) = (x_1 * x_2) \perp (x_1 * x_3).$$

Using (6.7) and the part (i) of our proposition, we obtain

$$(6.9) \qquad \begin{aligned} &x_1 * (x_2 \perp x_3) = \alpha x_1 ((x_2 \perp x_3)) + \beta ((x_2 \perp x_3)) x_1 = \alpha x_1 x_2 + \\ &\alpha x_1 x_3 + \beta x_2 x_1 + \beta x_3 x_1 + \alpha \gamma x_1 g(x_2, x_3) + \beta \gamma g(x_2, x_3) x_1. \end{aligned}$$

On the other hand, we get

$$(6.10) \qquad \begin{aligned} &(x_1 * x_2) \perp (x_1 * x_3) = (\alpha x_1 x_2 + \beta x_2 x_1) \perp (\alpha x_1 x_3 + \beta x_3 x_1) = \alpha x_1 x_2 + \\ &\alpha x_1 x_3 + \beta x_2 x_1 + \beta x_3 x_1 + \gamma g(\alpha x_1 x_2 + \beta x_2 x_1, \alpha x_1 x_3 + \beta x_3 x_1). \end{aligned}$$

Comparing the degrees in the expressions (6.9) and (6.10) we obtain

$$(6.11) \qquad deg\,(\alpha \gamma x_1 g(x_2, x_3) + \beta \gamma g(x_2, x_3) x_1) = k + 1$$

and

$$(6.12) \qquad deg\,\gamma \cdot g(\alpha x_1 x_2 + \beta x_2 x_1, \alpha x_1 x_3 + \beta x_3 x_1) > k + 1.$$

This contradiction shows that $\gamma = 0$ and (ii) follows.

(iii) By definition of operation $\circ$ in $F^*$, we have

(6.13) $$\alpha \circ x_1 = s(\alpha x_1), \ \alpha \in K, \ x_1 \in X$$

Since, according to (6.2), $sF(x_1) = F(x_1)$, and so

(6.14) $$\alpha \circ x_1 = P_\alpha(x_1),$$

where $P_\alpha(x_1) \in F(x_1)$. Consider the equality

(6.15) $$\alpha \circ (x_1 + x_2) = \alpha \circ x_1 + \alpha \circ x_2.$$

By (6.13) we get

$$\alpha \circ (x_1 + x_2) = \theta_{x_1 + x_2} s(\alpha x_1) = \theta_{x_1 + x_2} P_\alpha(x_1) = P_\alpha(x_1 + x_2).$$

From (6.15) follows

(6.16) $$P_\alpha(x_1 + x_2) = P_\alpha(x_1) + P_\alpha(x_2).$$

Write $P_\alpha(x_1)$ as an element of $U_1 = U(F(x_1))$: $P_\alpha(x_1) = \sum_i \alpha_i x_1^{k_i}, \ \alpha_i \in K, \ k_i > 0$.
By (6.16) we get

(6.17) $$\sum_i \alpha_i (x_1 + x_2)^{k_i} = \sum_i \alpha_i x_1^{k_i} + \sum_i \alpha_i x_2^{k_i}$$

Since $U(F(x_1, x_2))$ is a graded algebra, $(x_1 + x_2)^{k_i} = x_1^{k_i} + x_2^{k_i}$ for every $k_i$. Hence all $k_i = 1$. Thus, $P_\alpha(x_1) = s(\alpha x_1) = \alpha_1 x_1$ for some $\alpha_1 \in K$. Denote $\varphi(\alpha) = \alpha_1$, where $\varphi : K \to K$ is a mapping of $K$ into $K$. Since $s : F(x_1) \to F(x_1)$ is a bijection on $F(x_1)$, the mapping $\varphi$ is a bijection on $K$. From $\xi \circ (\mu \circ x_1) = (\xi \mu) \circ x_1$ and $(\xi + \mu) \circ x_1 = \xi \circ x_1 + \mu \circ x_1$, $\xi, \mu \in K$, follows $P_{\xi\mu}(x_1) = P_\xi P_\mu(x_1)$ and $P_{\xi+\mu}(x_1) = P_\xi(x_1) + P_\mu(x_1)$. Thus $\varphi$ is an automorphism of $K$ as required. This ends the proof. □

Now we consider a similar assertion for $\mathcal{A}_m$-varieties of algebras.

**Proposition 6.7.** Let $\Psi$ be an $\mathcal{A}_m$-variety of algebras over a ring $K$, $F_m$ be a finitely generated free algebra of $\Psi$, $\Phi$ be a quasi-inner automorphism of a semigroup $End\, F_m$ and $F_m^*$ be the derivative algebra associated with $\Phi$. The following statements hold:

(a) If $|K| = \infty$ or $|K| = p^k, k \geq 1$, where $p^k \nmid m - 1$ and $m \neq 2$, then the conclusions (i), (ii), (iii) of Proposition 6.6 are fulfilled for $F_m^*$.

(b) If $|K| = p^k, k \geq 1$, where $p^k | m - 1$, then the conclusions (i), (iii) of Proposition 6.6 are fulfilled for $F_m^*$ but instead (ii) the following is true

(ii)′ For any $a, b \in F_m^*$ we have $a \perp b = a + b + \gamma g(a, b)$, $\gamma \in K$, such that if $\gamma \neq 0$ then $g(x_1, x_2)$ is a nonzero homogeneous symmetric polynomial of degree $m - 1$ in $U(F_m)$ satisfied the following system of functional equations:

(6.18) $$g(x_1, x_2) + g(x_1 + x_2, x_3) = g(x_2, x_3) + g(x_1, x_2 + x_3)$$
$$g(x_1, -x_1) = 0,$$

where $x_i \in X$, $i \in [13]$.

(c) If $K$ any domain and $m = 2$, then the conclusion (ii) and (iii) of Proposition 6.6 are fulfilled for $F_2^*$ but instead (i) holds: $a * b = 0$ for any $a, b \in F_2$.

*Proof.* We consider the natural embedding $F_m \subseteq U(F_m)$, where $U(F_m)$ is the associative free $m$-nilpotent algebra. Recall that $U(F_m)$ is isomorphic to $A(X)/I$, where $A(X)$ is a free associative algebra generated by $X = \{x_1, ..., x_n\}$ and $I$ an ideal generated as a verbal ideal by elements $x_1 \cdot ... \cdot x_m$. It is easy to check that Lemmas 6.2 and 6.5 remain valid for the algebras $F$ and $U(F_m)$. In the proofs of these Lemmas, it is sufficient to restrict our consideration to $m$-reduced polynomials in the algebra $U(F_m)$.

Let assumptions of the part (a) be fulfilled. In the same manner as in Proposition 6.6 we derive the same conclusions (i), (iii) from this Proposition. Let us prove the equality (ii): $a \perp b = a + b$, $a, b \in F_m^*$. We will use the notations from Proposition 6.6 and the operation $\perp$ as defined in (6.6). We have in $F_m^*$

$$(\alpha \circ x_1) \perp (\alpha \circ x_2) = \alpha \circ (x_1 \perp x_2), \ \alpha \in K, \ x_1, x_2 \in X.$$

Since the part (iii) from Proposition 6.6 is fulfilled in our case, we get

$$(6.19) \qquad (\varphi(\alpha) - \varphi(\alpha)^{m-1}) g(x_1, x_2) = 0$$

Thus the conditions of part (a), $|K| = \infty$ or $|K| = p^k, k \geq 1$, where $p^k \nmid m-1$, give us $g(x_1, x_2) = 0$. This proves (ii) and therefore, the part (a) of our Proposition.

Let assumptions of the part (b) be fulfilled. Again the conclusions (i), (iii) hold in this case and it is remains to check the conclusion (ii)$'$. Since $|K| = p^k, k \geq 1$, where $p^k \mid m-1$, it is possible that $g(x_1, x_2) \neq 0$. Let us analyze the proof of the part (ii) of Proposition 6.6 in this case. There we have obtained a contradiction comparing the degrees of polynomials in both parts of the equalities (6.10) and (6.9). As is easy to see there is no contradiction if and only if in these equalities either the polynomial $g(x_1, x_2)$ is homogeneous of degree $m-1$ or $g(x_1, x_2) = 0$.

Let $g(x_1, x_2)$ be a homogeneous polynomial of degree $m-1$. Since

$$(6.20) \qquad (x_1 \perp x_2) \perp x_3 = x_1 \perp (x_2 \perp x_3)$$

we have

$$G(G(x_1 x_2), x_3) = G(x_1, G(x_2, x_3));$$

$$G(x_1 + x_2 + \gamma g(x_1, x_2), x_3) = G(x_1, x_2 + x_3 + \gamma g(x_1, x_2));$$

$$\gamma g(x_1, x_2) + \gamma g(x_1 + x_2 + \gamma g(x_1, x_2), x_3) = \gamma g(x_2, x_3) + \gamma g(x_1, x_2 + x_3 + \gamma g(x_2, x_3)).$$

Since $g(x_1, x_2)$ is a homogeneous polynomial of degree $m-1$, we get

$$(6.21) \qquad g(x_1, x_2) + g(x_1 + x_2, x_3) = g(x_2, x_3) + g(x_1, x_2 + x_3)$$

Since $x_1 \perp x_2 = x_2 \perp x_1$, we have $g(x_1, x_2) = g(x_2, x_1)$, i.e., the polynomial $g(x_1, x_2)$ is symmetric. Now, let $s$ be the bijection adjoint to $\Phi$. Since $s(0) = 0$, we have $x_1 \perp (-x_1) = 0$, i.e., $g(x_1, -x_1) = 0$. This proves the part (b) of our assertion.

Finally, let $m = 2$, i.e. $F = F_2$ be an algebra with trivial multiplication. Since $F_2 \simeq F_2^*$, we have $a * b = 0$, $a, b \in F_2^*$. From $g(x_1, x_2) = 0$ follows $a \perp b = a + b$. Analogously, we may check that $\xi \circ a = \varphi(\xi) a$ for any $a \in F_2^*$, $\xi \in K$ and an automorphism $\varphi : K \to K$. This ends the proof.   $\square$

Now we will describe the case where the system of the functional equations (6.18) has a non-trivial solution.

**Lemma 6.8.** *The system (6.18) has a non-trivial solution in the class of homogeneous symmetric polynomials from $U(F_m)$ of degree $m-1$ over a ring $K$ iff*

$m = 2k$, $k \geq 1$, or char $K = 2$ and $m$ is arbitrary. In these two cases the polynomial

$$(6.22) \qquad g(x_1, x_2) = (x_1 + x_2)^{m-1} - x_1^{m-1} - x_2^{m-1}$$

is a non-zero particular solution of the system (6.18).

*Proof.* Let $g = g(x_1, x_2)$ be a solution of the system (6.18). Setting $x_2 = -x_3$ in the first equation of (6.18) and taking into account the symmetry of polynomial $g(x_1, x_2)$, we obtain

$$(6.23) \qquad g(x_1, x_2) = -g(x_1 + x_2, -x_1)$$

Consider a linear operator $\mathcal{S}$ on the algebra $U(F_m)$:

$$\mathcal{S}(f)(x_1, x_2) = -f(\Lambda(x_1, x_2)), \ \forall f \in U(F_m)),$$

where $\Lambda = \begin{pmatrix} 1 & 1 \\ 0 & -1 \end{pmatrix}$. By (6.23), $\mathcal{S}(g) = g$. Since $\Lambda^3 = -I$, we have

$$g(x_1, x_2) = \mathcal{S}^3(g)(x_1, x_2) = -g(-x_1, -x_2).$$

Thus, $g(x_1, x_2) = (-1)^m g(x_1, x_2)$. Therefore, if $m = 2k + 1$, $k \geq 0$ and char $K \neq 2$, we obtain $g(x_1, x_2) = 0$. Otherwise, a straightforward check shows that the polynomial (6.22) is a non-trivial solution of the system (6.18). □

## 7. Proofs of the main theorems

*Proof of Theorem A.* Let $\Phi \in Aut\, End\, F(X)$ be an automorphism of $End\, F$. By Proposition 5.3, $\Phi$ is a quasi-inner automorphism of $End\, F$. Let $s : F \to F$ be the bijection adjoint to an automorphism $\Phi$ and $F^*$ be the derivative algebra associated with $\Phi$. Recall that under our assumptions $s(x_i) = x_i$ for all $x_i \in X$.

Denote by $[\,,\,]_1$ and $[\,,\,]_2$ the Lie operations in the Lie algebras $F$ and $F^*$, respectively, i.e.,

$$F = \langle F; \ \cdot, +, [\,,\,]_1, 0 \rangle \ \text{ and } \ F^* = \langle F; \ \circ, \perp, [\,,\,]_2, 0 \rangle$$

By PBW Theorem the variety of Lie algebras is an $\mathcal{A}$-variety [1, 7]. Using the Proposition 6.6, we can write

$$[a, b]_2 = \alpha a \cdot b + \beta b \cdot a, \ \forall a, b \in F, \ \forall \alpha, \beta \in K,$$

where the operation $\cdot$ is a multiplication in the universal enveloping algebra $U(F)$. Since $[a, a]_2 = 0$, we have $\alpha = -\beta$, i.e., $[a, b]_2 = \alpha[a, b]_1$. Since $F^*$ is a free Lie algebra, $\alpha \neq 0$. By virtue of Proposition 6.1,

$$(7.1) \qquad s[a, b]_1 = [s(a), s(b)]_2 = \alpha[s(a), s(b)]_1, \ \forall a, b \in F$$

In the same way we obtain

$$(7.2) \qquad \begin{array}{c} \forall a, b \in F, \ \forall \xi \in K \\ s(a + b) = s(a) \perp s(b) = s(a) + s(b), \ \ s(\xi a) = \xi \circ s(a) = \varphi(\xi)s(a). \end{array}$$

Let $\xi_\alpha : F \to F$ be a bijection on $F$ defined by $\xi_\alpha(a) = \alpha a$ for any $a \in F$. Denote by $\tilde{s} = \xi_{\alpha^{-1}} s$ a bijection on $F$. It is evident that $\tilde{s}$ is a semi-inner automorphism of $End\, F$. Finally, we have

$$\Phi(\nu) = s\nu s^{-1} = \tilde{s}\nu \tilde{s}^{-1}, \ \forall \nu \in End\, F,$$

where $\tilde{s}$ is a semi-inner automorphism. This end the proof. □

*Remark* 7.1. Let $H_1$ and $H_2$ be Lie algebras over an $R_1$MF-domain $K$. Using the same arguments as in [24] and Theorem A one can prove that the following conditions are equivalent:

1. The algebras $H_1$ and $H_2$ are categorically equivalent.

2. The algebras $H_1^\sigma$ and $H_2$ are geometrically equivalent for some $\sigma \in Aut\,K$.

Here the algebras $H_1^\sigma$ and $H_2$ coincide as rings and multiplication by a scalar in $H_1^\sigma$ is defined by the rule:

$$\lambda \circ a = \lambda^{\sigma^{-1}} \cdot a,\ \forall \lambda \in K,\ \forall a \in H_1^\sigma.$$

Earlier, for the variety of Lie algebras over infinite fields this result was obtained in [24].

*Proof of Theorem B.* 1. Let us prove the part 1 of Theorem B. Since $\Phi$ is an automorphism of $End\,F_m$, by Proposition 5.3, $\Phi$ is quasi-inner. As above, we may consider a derivative algebra $F_m^* = \langle F_m;\ *, \circ, \perp, 0 \rangle$ associated with $\Phi$. By Proposition 6.7, we have in $F_m^*$

$$(7.3) \qquad \exists \alpha, \exists \beta \in K,\ \forall a, b \in F_m,\ a*b = \alpha a \cdot b + \beta b \cdot a,$$

Since $m > 2$, we have $\alpha^2 + \beta^2 \neq 0$. Taking into account the law of associativity in $F_m$: $(a*b)*c = a*(b*c)$, we arrive at $\alpha\beta = 0$. Since $K$ is a domain, $\alpha = 0$ or $\beta = 0$. Since $\alpha^2 + \beta^2 \neq 0$, we have $a*b = \alpha a \cdot b$, $\alpha \neq 0$, or $a*b = \beta b \cdot a$, $\beta \neq 0$.

By Proposition 6.1, the adjoint bijection $s$ to $\Phi$ is an isomorphism of $F_m$ into $F_m^*$. Since $|K| = p^k$, $k \geq 1$, $p \neq 2$, and $p^k|m-1$, by Proposition 6.7, part (b), we have for the bijection $s: F_m \to F_m$ the following equalities:

$$(7.4) \qquad \begin{array}{l} \exists \gamma \in K, \exists \alpha \neq 0 \in K, \exists \varphi \in Aut\,K, \forall \xi \in K, \forall a \forall b \in F_m, \\ s(a+b) = s(a) \perp s(b) = s(a) + s(b) + \gamma g(s(a), s(b)), \\ s(\xi a) = \xi \circ s(a) = \varphi(\xi) s(a),\ s(a \cdot b) = s(a) * s(b) = \alpha s(a) \cdot s(b), \\ (\text{or } \exists \beta \neq 0 \in K, \forall a, b \in F_m,\ s(a \cdot b) = s(a) * s(b) = \beta s(b) \cdot s(a)), \end{array}$$

such that if $\gamma \neq 0$ then $g(x_1, x_2)$ is a nonzero homogeneous symmetric polynomial of degree $m-1$ in $U(F_m)$ which satisfies the system (6.18).

Let $m \neq 2k$, $k \geq 1$, and $char\,K \neq 2$. By Lemma 6.8, we get $g(x_1, x_2) = 0$. If $\alpha \neq 0$ and $\beta = 0$ in (7.4), we obtain, as above in Theorem A, that the automorphism $\Phi$ is semi-inner. If $\alpha = 0$ and $\beta \neq 0$, the automorphism $\Phi$ is a composition of mirror and semi-inner automorphisms. Therefore, the group $Aut\,End\,F$ is generated by semi-inner and mirror automorphisms.

Let $m = 2k$, $k \geq 1$, or $char\,K = 2$. By Lemma 6.8, there exists a non-zero homogeneous symmetric polynomial $g = g(x_1, x_2)$ of degree $m-1$ in $U(F_m)$ which is a solution of (6.18). Now we take $g = g(x_1, x_2)$ and set

$$(7.5) \qquad \forall a \forall b \in F_m, \exists \alpha \neq 0 \in K,\ s(a \cdot b) = \alpha s(a) \cdot s(b)$$

in (7.4). Since

$$s(a+b) = s(a) + s(b) + \gamma g(s(a), s(b)),$$

we obtain

$$s^r(a+b) = s^r(a) + s^r(b) + r\gamma g(s^r(a), s^r(b)),\ r \geq 1.$$

From the last equality follows

$$(7.6) \qquad s^p(a+b) = s^p(a) + s^p(b),$$

whereas

$$s^{p-1}(a+b) \neq s^{p-1}(a) + s^{p-1}(b).$$

We obtain, as in Theorem A, that $\Phi^p$ is a semi-inner automorphism, whereas $\Phi^{p-1}$ is not. Hence, the automorphism $\Phi$ is $p$-semi-inner. Setting in (7.4)

$$\forall a \forall b \in F_m, \exists \beta \neq 0 \in K, \ s(a \cdot b) = \beta s(b) \cdot s(a),$$

and choosing $\varphi$ the identical automorphism of $K$, in the same way as above we come to a $p$-mirror automorphism of $End\,F_m$.

Now we have to prove the existence of $p$-semi-inner and $p$-mirror automorphisms of $End\,F_m$. To this end we define a mapping $\tau : F_m \to F_m$ such that

(7.7)
$$\exists \varphi \in Aut\,K, \forall \xi \in K, \forall a \forall b \in F_m, \forall x_i \in X,$$
$$\tau(0) = 0, \ \tau(x_i) = x_i, \ \tau(x_{i_1}^{k_1}...x_{i_r}^{k_r}) = x_{i_1}^{k_1}...x_{i_r}^{k_r},$$
$$\tau(\xi a) = \varphi(\xi)\tau(a), \tau(a+b) = \tau(a) + \tau(b) + g(\tau(a), \tau(b)),$$

where the polynomial $g = g(x_1, x_2)$ is a non-zero homogeneous symmetric polynomials of degree $m - 1$ in $U(F_m)$ which is a solution of (6.18). Since

$$g(a, b) + g(a + b, c) = g(b, c) + g(a, b + c), \ \forall a \forall b \forall c \in F_m,$$

we obtain

$$\tau((a + b) + c) = \tau(a + (b + c)).$$

Hence, the mapping $\tau$ is defined correctly. As above (see (7.6)), we obtain

(7.8)
$$\tau^p(a + b) = \tau^p(a) + \tau^p(b),$$

whereas

$$\tau^{p-1}(a + b) \neq \tau^{p-1}(a) + \tau^{p-1}(b).$$

From (7.8) and the condition $\tau(x_i) = x_i, \in [1n]$ follows $\tau^p = Id_{F_m}$, where $Id_{F_m}$ is the identical mapping on $F_m$. Therefore, $\tau$ is a bijection on $F_m$. Thus, there exists a bijection $\tau^{-1} : F_m \to F_m$ and, furthermore, it is easy to check that

(7.9)   $$\tau^{-1}(a + b) = \tau^{-1}(a) + \tau^{-1}(b) + (p-1)g(\tau^{-1}(a), \tau^{-1}(b)), \forall a \forall b \in F_m$$

Define a mapping $\Psi$ of $End\,F_m$ such that $\Psi(\nu) = \tau^{-1}\nu\tau, \forall \nu \in End\,F_m$. Using (7.9), we obtain that $\Psi$ is an automorphism of $End\,F_m$. Now it is clear that $\Psi^p$ is a semi-inner automorphism, whereas $\Psi^{p-1}$ is not. Hence, the automorphism $\Psi$ is $p$-semi-inner.

Setting

$$\tau(x_{i_1}^{k_1}...x_{i_r}^{k_r}) = x_{i_r}^{k_r}...x_{i_1}^{k_1}, \ x_{i_k} \in X$$

in the definition (7.7) and choosing $\varphi$ the identical automorphism of $K$, we come to a $p$-mirror automorphism of $End\,F_m$.

Now it is clear that the group $Aut\,End\,F_m$ is generated by semi-inner, mirror, $p$-semi-inner and $p$-mirror automorphisms. This proves the part 1 of Theorem B.

2. Let us prove the part (2) of Theorem B. Since $|K| = \infty$ or $|K| = p^k, k \geq 1$, where $p^k \nmid m - 1$, by Proposition 6.7, part (a), we have the following equalities for the bijection $s : F_m \to F_m$:

$$\exists \alpha \neq 0 \in K, \exists \varphi \in Aut\,K, \forall \xi \in K, \forall a \forall b \in F_m$$
$$s(a + b) = s(a) + s(b), \ s(\xi a) = \varphi(\xi)s(a), \ s(a * b) = \alpha s(a) \cdot s(b),$$
$$(or\ \exists \beta \neq 0 \in K, \forall a, b \in F_m, \ s(a \cdot b) = s(a) * s(b) = \beta s(b) \cdot s(a)).$$

In the same manner as above we can prove that every automorphism $\Phi$ of the semi-group $End\,F_m$ is either a semi-inner or a mirror automorphism, or a composition of them. This proves the part 2 of Theorem B.

3. In the case of the variety $\mathcal{N}_2$, the multiplication in algebra $F_2$ is trivial, i.e., $a \cdot b = 0$. Thus, we have the following equalities for the bijection $s : F_2 \to F_2$:

$$\exists \varphi \in Aut\, K, \forall \xi \in K, \forall a \forall b \in F_m$$
$$s(a + b) = s(a) + s(b), \; s(\xi a) = \varphi(\xi) s(a).$$

As above we obtain that all automorphisms of $End\, F_2$ are semi-inner. The proof is complete. $\qquad\square$

*Proof of Theorem C.* Let $\varphi \in Aut\, \mathcal{L}^\circ$. It is clear that $\varphi$ is an equinumerous automorphism. By Proposition 2.8, $\varphi$ can be represented as the composition of a stable automorphism $\varphi_S$ and an inner automorphism $\varphi_I$. Since a stable automorphism does not change free algebras from $\Theta^0$, we obtain that $\varphi_S \in Aut\, End\, F(x_1, ..., x_n)$, where $F = F(x_1, ..., x_n)$ is a finitely generated free Lie algebra of $\mathcal{L}$. By Theorem A, $\varphi_S$ is a semi-inner automorphism of $End\, F(x_1, ..., x_n)$. Using this fact and Reduction Theorem [15, 21] we obtain that the composition $\varphi = \varphi_S \varphi_I$ is a semi-inner automorphism of $\Theta^0$. This ends the proof. $\qquad\square$

Now we provide an example of the variety $\mathcal{N}_3$ of 3-nilpotent associative algebras over the field $\mathbf{F}_2$ with a free algebra $F_3 = F_3(x_1, ..., x_n)$ such that the group $Aut End\, F_3$ contains 2-inner and 2-mirror automorphisms.

**Example 7.2.** Consider the polynomial $g(x_1, x_2) = x_1 x_2 + x_2 x_1$ in $F_3$. The polynomial $g(x_1, x_2)$ is a solution of the system (6.18) in the class of homogeneous symmetric polynomials from $F_3$ of degree 2 (see also Lemma 6.8). In a similar way as in proof of Theorem B we can construct a bijection $s : F \to F$ with the help of the polynomial $g(x_1, x_2)$ so that

(7.10)
$$\forall \xi \in K, \forall a, b \in F_m, \; \forall x_i \in X,$$
$$s(0) = 0, \; s(x_i) = x_i, \; s(x_{i_1}^{k_1} ... x_{i_r}^{k_r}) = x_{i_1}^{k_1} ... x_{i_r}^{k_r},$$
$$s(a + b) = s(a) + s(b) + g(s(a), s(b)) = s(a) + s(b) + s(a)s(b) + s(b)s(a).$$

As in Theorem B, we can prove correctness of definition of the mapping $s$. It is clear that $s^2$ is an automorphism of $F_3$, whereas $s$ is not. Thus, the automorphism $\Phi$ of the semigroup $End\, F_3$, such that $\Phi(\nu) = s\nu s^{-1}$ for any $\nu \in End\, F_3$, is 2-inner. Setting in (7.10)

$$s(x_{i_1}^{k_1} ... x_{i_r}^{k_r}) = x_{i_r}^{k_r} ... x_{i_1}^{k_1}, \; x_{i_k} \in X,$$

we arrive at a 2-mirror automorphism of $End\, F_3$.

## 8. Quasi-inner automorphisms of the semigroup $End\, A(X)$

Let $A = A(x_1, ..., x_n)$ be a free associative algebra over a ring $K$. A description of $Aut\, End\, A(x_1, x_2)$, where $A(x_1, x_2)$ is a free two generated associative algebra over an infinite field, has been obtained in [4]. In this connection the following assertion presents interest

**Proposition 8.1.** Let $\Phi \in Aut\, End\, A$ be a quasi-inner automorphism of $End\, A$, where $A = A(x_1, ..., x_n)$, $n \geq 2$, be a free finitely associative algebra over a domain $K$. Then $\Phi$ is either a semi-inner or a mirror automorphism, or a composition of them.

*Proof.* Let $s : A \to A$ be the adjoint bijection to $\Phi$. By Proposition 6.6 we have for $s$ the following equalities

$$\exists \alpha \neq 0 \in K, \exists \varphi \in Aut\, K, \forall \xi \in K, \forall a \forall b \in A$$
$$s(a+b) = s(a) + s(b), \ s(\xi a) = \varphi(\xi) s(a), \ s(a*b) = \alpha s(a) \cdot s(b),$$
$$(or\ \exists \beta \neq 0 \in K, \forall a, b \in A, \ s(a \cdot b) = s(a) * s(b) = \beta s(b) \cdot s(a)).$$

In the same manner as in Theorem B, part 2, we can prove that $\Phi$ is either a semi-inner or a mirror automorphism, or a composition of them. $\qquad\square$

## 9. Acknowledgments

The author is grateful to his teacher B. Plotkin for most helpful discussions and comments on this work. The author thanks his colleagues G. Belitskii, A. Belov-Kanel, G. Mashevitzky, E. Plotkin and G. Zhitomirski for helpful comments. In particular, Example 4.5 and Lemma 6.8 are due to G. Belitskii.

## References

[1] Yu. Bahturin, *Identical relations in Lie algebras*, (VNU Science Press, Utrecht, 1987).

[2] A. Berzins, B. Plotkin, E. Plotkin, Algebric geometry in varieties of algebras with the given algebra of constants, *Journal of Math. Sciences* **102** (3) (2000) 4039-4070.

[3] A. Berzins, The group of automorphisms of the category of free associative algebra, *Preprint* (2004).

[4] A. Berzins, The group of automorphisms of semigroup of endomorphisms of free commutative and free associative algebra, *Preprint*, (2004)

[5] G. Birkhoff, Representability of Lie algebras and Lie groups by matrices, *Ann. of Math.* **38** (2) (1937) 526-532

[6] P. Cohn, *Free rings and their relations*, (Academic Press, London, 1985).

[7] V. Drensky, *Free algebras and PI-algebras*, (Springer-Verlag, Singapore, 2000).

[8] E. Formanek, A question of B. Plotkin about the semigroup of endomorpjsms of a free group, *Proc. American Math. Soc.* **130** (2001) 935-937.

[9] L. Gluskin, Automorphisms of multiplicative semigroup of matrix algebras (Russian), *Uspehi Mat. nauk* **67**(11) (1954) 42-56.

[10] F. Halezov, Automorphisms of matrix semigroups (Russian), *Ivanov Ped. Inst., Zap. Fiz.-Mat. Nauki* **5** (11) (1956) 199-206).

[11] Jose A. Hermida-Alonso, On linear algebra over commutative rings, in *Handbook of algebra* Vol. 3 (North-Holland, Amsterdam, 2003), pp. 3-61.

[12] I. Isaacs, Automorphisms of matrix algebras over commutative rings, *Linear Algebra Appl.* **31** (1980) 215-231.

[13] M. Jodiet, T. Y. Lam, Multiplicative maps of matrix semigroups, *Arch. Math.* **20** (1969) 10-16.

[14] I. Kaplansky, Elementary divisors and modules, *Trans. Amer. Math.* **66** (1949) 464-491.

[15] Y. Katsov, R. Lipyanski, B. Plotkin, Automorphisms of categories of free modules and free Lie algebras , (2004) pp. 18, to appear.

[16] R. Lipyanski, B. Plotkin, Automorphisms of categories of free modules and free Lie algebras, *Preprint. Arxiv:math. RA//0502212* (2005).

[17] S. Mac Lane, *Categories for the Working Mathematician*, (New York-Berlin: Spinger-Verlag, 1971).

[18] A. Malcev, On algebras defined by identities (Russian), *Mat. Sb.* **26** (1950) 19-33.

[19] G. Mashevitzky, Automorphisms of the semigroup of endomorphisms of free ring and free associative algebras, Preprint

[20] G. Mashevitzky, B. Schein, Automorphisms of the endomorphism semigroup of a free monoid or a free semigroup, *Proc. Amer. math. Soc.* **8** (2002) 1-10.

[21] G. Mashevitzky, B. Plotkin, E. Plotkin, Automorphisms of the category of free Lie algebras, *Journal of Algebra* **282** (2004) 490-512.

[22] G. Mashevitzky, B. Plotkin, E. Plotkin, Automorphisms of the category of free algebras of varieties, *Electron. Res. Announs. Amer. Math. Soc.* **8** (2002) 1-10.

[23] B. Plotkin, Seven lectures in universal algebraic geometry, *Preprint. Arxiv:math. RA/0502212* (2002).

[24] B. Plotkin, Algebra with the same (algebraic geometry), in *Proc. of the Steklov Institut of Mathematics* **242** (2003) 176-207.

[25] B. Plotkin, G. Zhitomirskii, On automorphisms of categories of free algebras of some varieties. *Preprint. Arxiv:math. RA/0501331* (2005).

[26] K. Zhevlakov, A. Slinko, I. Shestakov, A. Sirshov, *Rings that are nearly associative (Russian)*, (Nauka, Moscow, 1978).

[27] G. Zhitomirskii, Automorphisms of the semigroup of all endomorphisms of free algebras, *Preprint. Arxiv:math. GM/0510230 v1* (2005).

*Department of Mathematics, Ben-Gurion University of the Negev, Beer Sheva, 84105, Israel*
*E-mail address*: `lipyansk@cs.bgu.ac.il`